\documentclass{amsart}
\usepackage{amsmath}
\usepackage{amsfonts}
\usepackage{amssymb}
\usepackage{amsthm}
\usepackage{amsbsy}
\usepackage{amscd}
\usepackage{latexsym}
\usepackage{mathtools}
\usepackage{wasysym}
\usepackage{bm}
\usepackage[all]{xy}
\usepackage{url}
\usepackage{color}
\usepackage[greek,english]{babel}
\usepackage{cmap}
\usepackage{ucs}
\usepackage[utf8x]{inputenc}
\usepackage{graphicx}
\usepackage{diagmac2} 
\usepackage{comment} 
\usepackage{scalerel}
\usepackage[euler]{textgreek}

\newtheorem{defn}{Definition}[section]

\newtheorem{lem}{Lemma}[section]
\newtheorem{prop}{Proposition}[section]
\newtheorem{cor}{Corollary}[prop]

\newcommand{\set}[2]{\{#1 ~|~#2\}}

\newcommand{\Lan}{\mathcal{L}}
\newcommand{\LanP}{\mathcal{L}^{+}}

\newcommand{\Var}{\mathcal{V}}
\newcommand{\VarL}{\mathcal{V}_{\mathcal{L}}}
\newcommand{\theory}{\Sigma_{\mathcal{S}}}

\newcommand{\ConS}[1]{\textbf{Cn}_{\mathcal{S}}(#1)}

\newcommand{\Forms}{\textbf{Fm}}
\newcommand{\FormsL}{\textbf{Fm}_{\mathcal{L}}}
\newcommand{\FormAl}{\mathfrak{F}_{\mathcal{L}}}
\newcommand{\alg}[1]{\textbf{#1}}
\newcommand{\mat}[1]{\textbf{#1}}
\newcommand{\LinS}{\textbf{Lin}_{\mathcal{S}}}

\newcommand{\cM}{\mathcal{M}}
\newcommand{\fS}[1]{\mathfrak{S}(#1)}

\newcommand{\Mod}[1]{\textbf{Mod}(#1)}

\newcommand{\vdashS}{\vdash_{\mathcal{S}}}
\newcommand{\card}[1]{\textit{card}(#1)}
\newcommand{\ext}{\prec_{p}}

\newcommand{\define}{\stackrel{\text{df}}{\Longleftrightarrow}}

\newcommand{\aLog}{\mathcal{S}}

\begin{document}

\title{On Matrix Consequence (Extended Abstract)}

\author{
Alexei Muravitsky
}
\maketitle
\begin{center}
 Louisiana Scholars' College \\
Northwestern State University\\
Natchitoches, LA 71497, U.S.A. \\
\emph{Email}: alexeim@nsula.edu	
\end{center}

\begin{abstract}
These results are a contribution to the model theory of matrix consequence. We give a semantic characterization of uniform and couniform consequence relations. The first concept was introduced by {\L}o\'{s} and Suszko (1958), and the second by W\'{o}jcicki (1970) to characterize structural consequence relations having single adequate matrices. In different variants of the single-adequate-matrix theorem of {\L}o\'{s} and Suszko (corrected by W\'{o}jcicki), the uniformity of consequence is used in conjunction with other properties (e.g. finitariness), and in the W\'{o}jcicki theorem uniformity and couniformity of consequence are used together. These properties have never been treated individually, at least in a semantic manner. We consider these notions from a purely semantic point of view and separately, introducing the notion of a uniform bundle/atlas and that of a couniform class of logical matrices. Then, we show that any uniform bundle defines a uniform consequence; and if a structural consequence is uniform, then its Lindenbaum atlas is uniform. Thus, any structural consequence is uniform if, and only if, it is determined by a uniform bundle/atlas. On the other hand, any couniform set of matrices defines a couniform structural consequence.  Also, the Lindenbaum atlas of a couniform structural consequence is couniform. Thus, any structural consequence is couniform if, and only if, it is determined by a couniform bundle/atlas. We then apply these observations to compare structural consequence relations that are defined in different languages when one language is a primitive extension of another. We obtain that for any structural consequence defined in a language having (at least) a denumerable set of sentential variables, if this consequence is uniform and couniform, then it and the \emph{ W\'{o}jcicki's consequence} corresponding to it, which is defined in any primitive extension of the given language, are determined by one and the same atlas which is both uniform and couniform.
\end{abstract}

\section{Preliminaries}
We consider a structural consequence relation, not necessarily finitary (or compact), in an abstract propositional language $\Lan$. We will be using a unifying term \textit{abstract logic} (usually denoted by $\aLog$) in order to employ in one and the same context the consequence relation $\vdashS$ and the consequence operator $\textbf{Cn}_{\mathcal{S}}$ both associated with $\aLog$.

The set of all $\Lan$-variables is denoted by $\VarL$, the set of all $\Lan$-formulae by $\FormsL$. For any $\alpha\in\FormsL$,
$\Var(\alpha)$ is the set of all variables occurring in $\alpha$. Further, if $\Var\subseteq\VarL$, we denote:
\[
\FormsL[\Var]:=\set{\alpha\in\FormsL}{\Var(\alpha)\subseteq\Var}.
\]
$\FormAl$ stands for a formula algebra in $\Lan$. Given a set $X\subseteq\FormsL$ and an abstract logic $\aLog$, $X$ is an $\aLog$-\textbf{\textit{theory}} if $\ConS{X}=X$.

It is well known (see, e.g.,~\cite{wojcicki1988}) that any structural consequence is determined by a class of (logical) matrices and vice versa. We remind that, given a logical matrix $\mat{M}=\langle\alg{A},D\rangle$, where $\alg{A}$ is an algebra of the same signature as the signature of $\Lan$ and $D\subseteq|\alg{A}|$ (we call such a matrix an $\Lan$-\textit{\textbf{matrix}} and $D$ its \textit{\textbf{logical filter}}), a \textit{\textbf{matrix consequence relative to}} \mat{M} (or \mat{M}-\textit{\textbf{consequence}} for short) is the following relation:
\[
X\models_{\mat{M}}\alpha~\define~(v[X]\subseteq D\Longrightarrow v[\alpha]\in D,~\text{for any valuation $v$ in \alg{A}}),
\]
where $X\cup\lbrace\alpha\rbrace\subseteq\FormsL$.

The last definition is extended to a (nonempty) class $\cM=\lbrace\mat{M}_i\rbrace_{i\in I}$ of $\Lan$-matrices as follows:
\begin{equation}\label{E:matrix-consequence}
X\vdash_{\mathcal{M}}\alpha~\define~(X\models_{\mat{M}_i}\alpha,~\textit{for each $i\in I$}).
\end{equation}
We call the last relation a \textit{\textbf{matrix consequence relative to}} $\cM$ (or $\cM$-\textit{\textbf{consequence}} for short).

Given an abstract logic $\aLog$, a logical matrix $\mat{M}$ is an $\aLog$-\textit{\textbf{model}} if
\[
X\vdashS\alpha~\Longrightarrow X\models_{\mat{M}}\alpha.
\]
The set of all $\aLog$-models is denoted by $\Mod{\aLog}$.\\

In fact, if a consequence relation is determined by a class of $\Lan$-matrices, it is also determined by a bundle or atlas.
\begin{defn}\label{D:bundle/atlas}
	A $($nonempty$)$ family {\em$\mathcal{B}=\lbrace\langle\alg{A}_{i},D_i\rangle\rbrace_{i\in I}$} of $\Lan$-matrices is called a \textbf{bundle} if for any $i,j\in I$, the algebras {\em$\alg{A}_i$} and {\em$\alg{A}_j$} are isomorphic. Ignoring the differences between the algebras {\em$\alg{A}_i$}, we write {\em$\mathcal{B}=\lbrace\langle\alg{A},D_i\rangle\rbrace_{i\in I}$}.  A pair {\em$\langle\alg{A},\lbrace D_i\rbrace_{i\in I}\rangle$}, where each $D_i$ is a logical filter in $($or of$)$ {\em$\alg{A}$}, is called an \textbf{atlas}. By definition, the \textbf{matrix consequence relative to an atlas}
	{\em$\langle\alg{A},\lbrace D_i\rbrace_{i\in I}\rangle$} is the same as the matrix consequence related the bundle {\em$\lbrace\langle\alg{A},D_i\rangle\rbrace_{i\in I}$}, which in turn is defined in the sense of~{\em \eqref{E:matrix-consequence}}. An atlas {\em$\LinS[\theory]=\langle\FormsL,\Sigma_{\mathcal{S}}\rangle$}, where $\Sigma_{\mathcal{S}}$ is the set of all $\aLog$-theories, is called a \textbf{Lindenbaum atlas} $($relative to $\aLog$$)$.
\end{defn}

The transition from a class $\cM$ of matrices to a bundle (or atlas), which would define the $\cM$-consequence, can be done (at least) in two ways.

An indirect way for such a transition, which is based on the fact that  $\cM$-consequence is structural, uses the Lindenbaum theorem; cf., e.g., W\'{o}jcicki~\cite{wojcicki1970}, theorem 2.4,~\cite{wojcicki1988} or~\cite{citkin-muravitsky2019}, corollary 3.3.5.1.

According to a direct method, given $\mathcal{M}=\lbrace\langle\alg{A}_i,D_i\rangle\rbrace_{i\in I}$, we define an atlas  $\mathcal{M}^{\ast}=\langle\alg{A},\lbrace D_{i}^{\ast}\rbrace_{i\in I}\rangle$ as follows:
\[
\begin{array}{rll}
\alg{A}:=\prod_{i\in I}\alg{A}_i
&\text{and, given $i\in I$,}~D_{i}^{\ast}:=\prod_{j\in I} H_j,
&\text{where}~H_j:=\begin{cases}
\begin{array}{cl}
D_i &\text{if $j=i$}\\
|\alg{A}_j| &\text{if $j\neq i$}.
\end{array}
\end{cases}
\end{array}
\]

We claim that
\[
X\vdash_{\mathcal{M}}\alpha~\Longleftrightarrow~X\vdash_{\mathcal{M}^{\ast}}\alpha.
\]

In~\cite{los-suszko1958}, the authors used the notion of uniform consequence (a \textit{\textbf{uniform abstract logic}} in our terminology) in an attempt to prove a criterion for a \textit{\textbf{single-matrix consequence}}, that is when a structural consequence can be determined by a single matrix. As W\'{o}jcicki (see, e.g.,~\cite{wojcicki1988}) has shown, in addition to the notion of uniformity, the notion of a \textit{\textbf{couniform abstract logic}} (in our terminology) is needed for such a criterion.

Originally, the concepts of uniformity and couniformity were defined syntactically. Below we treat them semantically.

\section{Uniform consequence}
For an arbitrary $\Lan$-matrix $\mat{M}=\langle\alg{A},D\rangle$, we define
an operator $\mathfrak{S}:\mat{M}\mapsto\fS{\mat{M}}$ where
\begin{equation}\label{E:S-operator-in-L}
\begin{array}{r}
X\in\mathfrak{S}(\mat{M})~\define~\text{there is an $\Lan$-valuation $v$ in \alg{A} such that}\\ 
X=\set{\alpha\in\FormsL}{v[\alpha]\in D}.
\end{array}
\end{equation}

We notice that
\[
X\models_{\mat{M}}\alpha~\Longleftrightarrow~\forall Y\in\fS{\mat{M}}.~X\subseteq Y~\Longrightarrow~\alpha\in Y.
\]

In case that we have a set $\cM=\lbrace\mat{M}_i\rbrace_{i\in I}$ of $\Lan$-matrices, the corresponding $\cM$-consequence can be reformulated in terms of $\mathfrak{S}$-operator as follows:
\begin{equation}\label{E:M-con-via-S-operator}
X\vdash_{\cM}\alpha~\Longleftrightarrow~\forall i\in I\,\forall Y\in\fS{\mat{M}_i}.~X\subseteq Y~\Longrightarrow~\alpha\in Y.
\end{equation}

We continue with the following definition.
\begin{defn}[uniform bundle/atlas]\label{D:uniform-bundle}
	Let {\em$\cM=\lbrace\mat{M}_i\rbrace_{i\in I}$} be a bundle.
	We say that $\cM$ is \textit{\textbf{uniform}} if for any {\em$X\cup Y\subseteq\FormsL$} with {\em$X\subseteq\FormsL[\VarL\setminus\Var(Y)]$}, the following condition is fulfilled:
	{\em\[
		\begin{array}{l}
		\forall i,j\in I\,\forall Z_i\in\fS{\mat{M}_i}\forall Z_j\in\fS{\mat{M}_j}\exists k\in I\exists Z_k\in\fS{\mat{M}_k}.\,X\subseteq Z_i\,\&\,Y\subseteq Z_j\subset\FormsL\\
		\Longrightarrow
		Z_i\cap\FormsL[\VarL\setminus\Var(Y)]=Z_k\cap\FormsL[\VarL\setminus\Var(Y)]\,\&\,
		Z_j\cap\FormsL[\Var(Y)]\subseteq Z_k.
		\end{array}
		\]}	
	An atlas is \textbf{uniform} if the corresponding bundle is uniform.
\end{defn}

The above definition is justified by the following two propositions.
\begin{prop}\label{P:M-uniform-con-is-uniform}
	Let {\em$\cM$} be a uniform bundle.
	Then the $\cM$-consequence is uniform.
\end{prop}
\begin{prop}\label{P:LinS-uniform}
	If a structural abstract logic $\aLog$ is uniform, then its Lindenbaum atlas {\em$\LinS[\theory]$} is uniform.
\end{prop}

\begin{cor}\label{P:uniformity-by-uniform-bundle}
	A structural abstract logic is uniform if, and only if, it is determined by a uniform bundle/atlas.
\end{cor}

\section{Couniform consequence}
We start with the following observation.
\begin{prop}\label{P:couniformity-criterion}
	Let $\aLog$ be a structural abstract logic that is determined by
	a set {\em$\cM=\lbrace\mat{M}_j\rbrace_{j\in J}$} of $\Lan$-matrices. The logic $\aLog$ is couniform if, and only if, for any collection $\lbrace X_i\rbrace_{i\in I}$ of formula sets with $\Var(X_i)\cap\Var(X_j)=\varnothing$, providing that $i\neq j$, and $\Var(\bigcup_{i\in I}\lbrace X_i\rbrace) \neq\Var_{\Lan}$, the following condition is satisfied:
	{\em\begin{equation}\label{E:couniformity-2}
		\begin{array}{r}
		\forall i\in I\exists j\in J\exists Z_j\in\mathfrak{S}(\mat{M}_j)(X_i\subseteq Z_j\subset\FormsL)~\Longrightarrow\\
		\exists k\in J\exists Z_k\in\mathfrak{S}(\mat{M}_k)(\bigcup_{i\in I}\lbrace X_i\rbrace\subseteq Z_k\subset\FormsL).
		\end{array}
		\end{equation}}	
\end{prop}

The last proposition leads to the following definition.
\begin{defn}\label{D:couniform-set-matrices} 
	A nonempty set {\em$\cM=\lbrace\mat{M}_j\rbrace_{j\in J}$} of $\Lan$-matrices is called \textbf{couniform} if for any collection $\lbrace X_i\rbrace_{i\in I}$ of formula sets with $\Var(X_i)\cap\Var(X_j)=\varnothing$, providing that $i\neq j$, and $\Var(\bigcup_{i\in I}\lbrace X_i\rbrace) \neq\Var_{\Lan}$, the condition~\eqref{E:couniformity-2} is satisfied. An atlas {\em$\langle\alg{A},\lbrace D_i\rbrace_{i\in I}\rangle$} is \textbf{couniform} if the corresponding bundle {\em$\lbrace\langle\alg{A},D_i\rangle\rbrace_{i\in I}$} is couniform.	
\end{defn}

Now we obtain
\begin{prop}\label{P:LinS-couniform}
	Let $\aLog$ be a structural couniform logic. Then its Lindenbaum atlas {\em$\LinS[\theory]$} is couniform.
\end{prop}

\begin{cor}
A structural logic is couniform if, and only if, it is determined by a couniform bundle/atlas.
\end{cor}

Using W\'{o}jcicki's criterion (see, e.g.,~\cite{wojcicki1988}, theorem 3.2.7, or~\cite{citkin-muravitsky2019}, proposition 4.1.2), we also obtain the following.
\begin{cor}
A structural abstract logic is determined by a single matrix if, and only if, it can be determined by a uniform bundle/atlas and a couniform bundle/atlas.
\end{cor}

\section{Extensions}
We start with the following two definitions.
\begin{defn}[extension and primitive extension of a language]
	Let $\Lan$ and $\Lan^{\prime}$ be two propositional languages with $\VarL\cup\mathcal{C}_{\mathcal{L}}\cup\mathcal{F}_{\mathcal{L}}\subseteq\Var_{\mathcal{L}^{\prime}}\cup\mathcal{C}_{\mathcal{L}^{\prime}}\cup\mathcal{F}_{\mathcal{L}^{\prime}}$, then $\Lan^{\prime}$ is called an \textbf{extension} of $\Lan$. If $\mathcal{C}_{\mathcal{L}}\cup\mathcal{F}_{\mathcal{L}}=\mathcal{C}_{\mathcal{L}^{\prime}}\cup\mathcal{F}_{\mathcal{L}^{\prime}}$ and $\VarL\subseteq\Var_{\mathcal{L}^{\prime}}$, then we say that $\Lan^{\prime}$ is a \textbf{primitive extension} of $\Lan$, symbolically {\em$\Lan\ext\LanP$}.	
\end{defn}

\begin{defn}
	Let $\aLog$ be an abstract logic in a language $\Lan$ and $\aLog^{\prime}$ be an abstract logic in a language $\Lan^{\prime}$ which is an extension of $\Lan$. Then $\aLog^{\prime}$ is called a \textbf{conservative extension of}
	$\aLog$, symbolically $\aLog\prec_{c}\aLog^{\prime}$, if for any set {\em$X\cup\lbrace\alpha\rbrace\subseteq\FormsL$},
	\[
	X\vdashS\alpha~\Longleftrightarrow~X\vdash_{\mathcal{S}^{\prime}}\alpha.
	\]
\end{defn}

\begin{prop}\label{P:equality-of-Mod-classes}
	Let $\aLog$ and $\aLog^{\prime}$ be structural abstract logics in languages $\Lan$ and $\Lan^{\prime}$, respectively, with $\Lan\prec_{p}\Lan^{\prime}$ and
	$\card{\Var_{\mathcal{L}^{\prime}}}=\card{\VarL}\ge\aleph_{0}$. If $\aLog\prec_{c}\aLog^{\prime}$, then {\em$\Mod{\aLog}=\Mod{\aLog^{\prime}}$}.
\end{prop}

\begin{defn}[consequence relation $\vdash_{(\mathcal{M})^{+}}$]
	Let $\mathcal{M}$ be a set of logical matrices of type $\Lan$. Assume that $\Lan\prec_{p}\LanP$. A consequence relation in $\LanP$ determined by the set $\mathcal{M}$, is denoted by $\vdash_{(\mathcal{M})^{+}}$.
\end{defn}

\begin{prop}
	Let $\mathcal{M}$ be a set of logical matrices of type $\Lan$. Assume that $\Lan\prec_{p}\LanP$. The consequence relation $\vdash_{(\mathcal{M})^{+}}$ if finitary if, and only if,
	$\vdash_{\mathcal{M}}$ is finitary.
\end{prop}

\begin{defn}[relation $\vdash_{(\mathcal{S})^{+}}$]\label{D:(S)^+}
	Let a language $\Lan^{+}$ be an extension of a language $\Lan$. Also, let $\mathcal{S}$ be an abstract logic in $\Lan$. We define a relation $\vdash_{(\mathcal{S})^{+}}$ in {\em$\mathcal{P}(\Forms_{\Lan^{+}})\times\Forms_{\Lan^{+}}$} as follows:
	{\em\[
		\begin{array}{rl}
		X\vdash_{(\mathcal{S})^{+}}\alpha\stackrel{\text{df}}{\Longleftrightarrow}
		\!\!\!&\textit{there is a set {\em$Y\cup\lbrace\beta\rbrace\subseteq\Forms_{\Lan}$} and an $\Lan^{+}$\!-substitution}\\
		&\textit{$\sigma$ such that $\sigma(Y)\subseteq X$, $\alpha=\sigma(\beta)$ and $Y\vdash_{\mathcal{S}}\beta$}.
		\end{array}
		\]}
\end{defn}

We borrow some results from~\cite{citkin-muravitsky2019}, collecting them in one lemma.
\begin{lem}[\cite{citkin-muravitsky2019}, Proposition3.5.4, Proposition 3.5.5, Lemma 4.1.2, Lemma 4.1.4]\label{L:lemma-in-Extensions}
	Let $\aLog$ be a structural logic which is both uniform and couniform in a language $\Lan$ with $\card{\Var_{\Lan}}\ge\aleph_{0}$ and $\Lan^{+}$ be a primitive extension of $\Lan$. Then $(\mathcal{S})^{+}$ is a structural abstract logic which is both uniform and couniform.	
\end{lem}

Based on Lemma~\ref{L:lemma-in-Extensions}, given an abstract logic $\aLog$,
we call $(\mathcal{S})^{+}$ a \textit{\textbf{W\'{o}jcicki logic}},
and the consequence related to the latter a \textit{\textbf{W\'{o}jcicki consequence}}.\\

Using the last lemma, we obtain the following.
\begin{prop}
	Let $\aLog$ be a structural logic which is both uniform and couniform in a language $\Lan$ with $\card{\Var_{\Lan}}\ge\aleph_{0}$ and $\Lan^{+}$ be a primitive extension of $\Lan$. Then $\aLog$ and $(\mathcal{S})^{+}$ are determined by one and the same atlas which is both uniform and couniform.
\end{prop}

\paragraph{\textbf{Acknowledgment}}
I am grateful to Alex Citkin for numerous discussions regarding the topic of this article.

%\bibliographystyle{plain}
%\bibliography{Bibtex}

\end{document}